\documentclass[12pt]{amsproc}
\usepackage{amsmath,amssymb}
\usepackage{pstricks}
\usepackage{pst-all}

\usepackage[T1]{fontenc}
\def\mynode#1{\Cnode[fillstyle=solid,fillcolor=black]{#1}}
\renewcommand{\a}{\alpha}
\renewcommand{\b}{\beta}

\title[Random walks and orthogonal polynomials]{Random walks and
orthogonal polynomials: some challenges}
\author{F. Alberto Gr\"unbaum}
\thanks{The author was supported in part by NSF Grant \# DMS 0603901}
\date{}
\address{Department of Mathematics \\ University of California \\ Berkeley,
CA\ \ 94720}

\begin{document}

\begin{abstract}
The study of several naturally arising ``nearest neighbours" random walks
benefits from the study of the associated orthogonal polynomials and their
orthogonality measure. I consider extensions of this approach to a larger class
of random walks. This raises a number of open problems.
\end{abstract}

\maketitle

\bigskip

To Henry, teacher and friend, with gratitude and admiration.

\section{Introduction}

Consider a birth and death process, i.e., a discrete time Markov
chain on the nonnegative integers, with a
one step transition probability matrix ${\mathbb P}$. There is then a
time-honored way of writing down the n-step transition
probability matrix ${\mathbb P}^n$ in terms of the orthogonal
polynomials associated to ${\mathbb P}$ and the spectral measure.
This goes
back to \cite{KMcG} and, as they observe, this is nothing but an
application of the spectral theorem.
One can find some precursors of these powerful ideas, see for
instance \cite{H,LR}.
In as much as this is such a deep and general result, it holds in
many setups, such as a nearest neighbours random walk on the $N$th-roots of unity.
In general this representation of ${\mathbb P}^n$
allows
one to relate properties of the Markov chain, such as recurrence or
other limiting behaviour, to properties of the orthogonality measure.

\bigskip
In the few cases when one can get one's
hands on the orthogonality measure and the polynomials themselves
this gives fairly explicit answers to various questions.

\bigskip
The two main drawbacks to the applicability of this representation
(to be recalled below) are:

\begin{itemize}
\item[a)] typically one cannot get either the polynomials or the
measure explicitly.
\item[b)] the method is restricted to ``nearest neighbour" transition
probability chains that give rise to tridiagonal matrices and thus
     to orthogonal polynomials.
\end{itemize}

\bigskip
The {\bf challenge} that we pose in this paper is very simple: try to extend
the class of Markov processes whose study can benefit from a similar association.

There is an important collection of papers that study in detail the cases where
the entries in ${\mathbb P}$ depends linearly, quadratically or even
rationally on the
index $n$. We make no attempt to review these results, but we just mention that
the linear case involves (associated) Laguerre and Meixner polynomials, and the
other cases involve associated dual Hahn polynomials. For a very
small sample of important sources
dealing with this connection see \cite{Ch,vD,ILMV}.

\bigskip
The plan for this paper is as follows:

\begin{itemize}
\item[1)] we first review briefly the approach of S.~Karlin and J.~McGregor.
This is done in section 2.
\item[2)] we consider a few examples of physically important Markov chains
that happen to feature rather well known families of orthogonal
polynomials. This is done in sections 3,4,5.
\item[3)] we propose a way of extending this representation to the case of
certain Markov chains where the one-step transition probabilty matrix
is not necessarily
tridiagonal. For concretness we restrict ourselves to the case of
pentadiagonal matrices or more generally block tridiagonal matrices.
This is illustrated with some examples. This material takes up sections 6,7,8,9 and 10.
\item[4)] A number of open problems are mentioned along the way. A few more
are listed in section 11.
\end{itemize}

\bigskip
Both parts 1) and 2) above are well known while the proposal in 3)
appears to be new.

\section{The Karlin-McGregor representation}

If we have
$${\mathbb P}_{i,j}= Pr\{X(n+1)=j | X(n)=i\}$$
for the 1-step transition probability of our Markov chain, and we
put $p_i={\mathbb P}_{i,i+1}$, $q_{i+1}={\mathbb P}_{i+1,i}$, and
$r_{i}={\mathbb P}_{i,i}$ we get
for the matrix ${\mathbb P}$, in the case of a birth and death
process, the expression
\[
{\mathbb P} = \begin{pmatrix}
r_0 & p_0 & 0 & 0 \\
q_1 & r_1 & p_1 & 0 \\
0 & q_2 & r_2 & p_2 \\
& &\ddots &\ddots &\ddots
\end{pmatrix}
\]
We will assume that $p_j > 0$, $q_{j+1} > 0$ and $r_j \ge 0$ for $j
\ge 0$.  We also assume $p_j + r_j + q_j = 1$ for $j \ge 1$ and by
putting $p_0 + r_0 \le 1$ we allow for the state $j=0$ to be an
absorbing state (with probability $1 - p_0 - r_0$). Some of these
conditions can be relaxed.

If one introduces the polynomials $Q_j(x)$ by the conditions
$Q_{-1}(0) = 0$, $Q_0(x) = 1$ and using the notation
\[
   Q(x) = \begin{pmatrix}
Q_0(x) \\
Q_1(x) \\
\vdots
\end{pmatrix}
\]
we insist on the recursion relation
\[
{\mathbb P} Q(x) = x Q(x)
\]
one proves the existence of a unique measure $\psi(dx)$ supported in
$[-1,1]$ such that
\[
\pi_j \int_{-1}^1 Q_i(x)Q_j(x)\psi(dx) = \delta_{ij}
\]
and one gets the Karlin-McGregor representation formula
\[
({\mathbb P}^n)_{ij} = \pi_j \int_{-1}^1 x^nQ_i(x)Q_j(x)\psi(dx).
\]
%where the constants $\pi_j$ are given by
%\[
%\pi_0 = 1,\quad \pi_j = \frac {p_0p_1 \dots p_{j-1}}{q_1q_2 \dots
%q_j} \quad j \ge 1.
%\]
Many general results can be obtained from the representation formula
given above.  Some of these will be given for certain examples in the next
three sections.

Here we just remark that the existence of
\[
\lim_{n \to \infty} ({\mathbb P}^n)_{ij}
\]
is equivalent to $\psi(dx)$ having no mass at $x=-1$.  If this is the
case this limit is positive exactly when $\psi(dx)$ has some mass at
$x = 1$.

If one notices that $Q_n(x)$ is nothing but the determinant of the
$(n+1) \times (n+1)$
upper-left corner of the matrix $x I -{\mathbb P}$, divided by the factor
$$p_0 p_1 ...p_{n-1},$$  and one defines the
polynomials $q_n(x)$ by
solving the same three-term recursion relation satisfied by the
polynomials $Q_n(x)$, but with the indices shifted by one, and
the initial conditions $q_0(x)=1$, $q_1(x)=(x-r_1)/p_1$ then it is
clear that the $(0,0)$ entry
of the matrix
\[
(xI -{\mathbb P})^{-1}
\]
should be given, except by the constant $p_0$, by the limit of the ratio
$$q_{n-1}(x)/Q_n(x).$$

On the other hand the same spectral theorem alluded to above
establishes an intimate relation between
$$lim_{n \to \infty} q_{n-1}(x)/Q_n(x)$$ and $$  \int_{-1}^1 \frac
{d\psi(\lambda)}{x- \lambda} .$$

We will see in some of the examples a probabilistic interpretation for the
right hand side of the expression above in terms of generating functions.

The same connection with orthogonal polynomials holds in the case of a
birth and death process with continuous time, and this has been extensively
described in the literature. The discrete time situation discussed
above is enough to
illustrate the power of this method.

\section{The Ehrenfest urn model}

Consider the case of a Markov chain in the finite state space ${0,1,2,...,2N}$
where the matrix ${\mathbb P}$ given by
\[
\begin{pmatrix}
0 & 1 \\
\frac {1}{2N} & 0 & \frac {2N-1}{2N} \\
& \frac {2}{2N} & 0 & \frac {2N-2}{2N} \\
& & \ddots & 0 & \ddots \\
& & & \ddots & \ddots & \ddots \\
& & & & \ddots  & 0 & \frac {1}{2N} \\
& & & & & \frac {2N}{2N} & 0
\end{pmatrix}
\]

This situation arises in a model introduced by P.\ and T.~Ehrenfest,
see \cite{E}, in an effort
to illustrate the issue that irreversibility and recurrence can coexist.
The background here is, of course, the famous $H$-theorem of L.~Boltzmann.

For a more detailed discussion of the
model see \cite{F,K}. This model has also been considered in dealing
with a quantum mechanical version of a
discrete harmonic oscillator by Schroedinger himself, see \cite{SchK}.

In this case the corresponding orthogonal polynomials (on a finite
set) can be given explicitly. Consider the so called
Krawtchouk polynomials, given by means of the (truncated) Gauss series
\[
{}_2{\tilde F}_1 \left( \begin{matrix}
a,b \\
c
\end{matrix}; z\right) = \sum_0^{2N} \frac {(a)_n(b)_n}{n!(c)_n} z^n
\]
with
\[
(a)_n\equiv a(a+1)\dots(a+n-1),\quad (a)_0 = 1
\]

The polynomials are given by
\begin{eqnarray*}
K_i(x) &= &{}_2{\tilde F}_1 \left( \begin{matrix}
-i,-x \\
-2N
\end{matrix}; 2\right) \\
x &= &0,1,\dots,2N;\quad i = 0,1,\dots,2N
\end{eqnarray*}

Observe that
\[
K_0(x) \equiv 1,K_i(2N) = (-1)^i.
\]
The orthogonality measure is read off from
\[
\sum_{x=0}^{2N} K_i(x)K_j(x) \frac {\begin{pmatrix}
2N \\
x
\end{pmatrix}}{2^{2N}} = \frac {(-1)^ii!}{(-2N)_i} \delta_{ij} \equiv
\pi_i^{-1} \delta_{ij} \quad 0 \le i,j \le 2N.
\]
These polynomials satisfy the second order difference equation
\[
\frac {1}{2} (2N-i)K_{i+1}(x) - \frac {1}{2} 2N K_i(x) + \frac {i}{2}
K_{i-1}(x) = -xK_i(x)
\]
and this has the consequence that
\[
\begin{pmatrix}
0 & 1 \\
\frac {1}{2N} & 0 & \frac {2N-1}{2N} \\
& \frac {2}{2N} & 0 & \frac {2N-2}{2N} \\
& & \ddots & 0 & \ddots \\
& & & \ddots & \ddots & \ddots \\
& & & & \ddots  & 0 & \frac {1}{2N} \\
& & & & & \frac {2N}{2N} & 0
\end{pmatrix} \begin{pmatrix}
K_0(x) \\
K_1(x) \\
{} \\
{} \\
{} \\
{} \\
{} \\
K_{2N}(x)
\end{pmatrix} = \left( 1 - \frac {x}{N}\right) \begin{pmatrix}
K_0(x) \\
K_1(x) \\
\vdots \\
K_{2N}(x)
\end{pmatrix}
\]
any time that $x$ is one of the values $x = 0,1,\dots,2N$.  This
means that the eigenvalues of the matrix ${\mathbb P}$ above are given by the
values of $1 - \frac {x}{N}$ at these values of $x$, i.e.,
\[
1,1 - \frac {1}{N},\dots,-1
\]
and that the corresponding eigenvectors are the values of
$[K_0(x),K_1(x),\dots,K_{2N}(x)]^T$ at these values of $x$.

Since the matrix ${\mathbb P}$ above is the one step transition
probability matrix for our urn model we conclude that
\[
({\mathbb P}^n)_{ij} = \pi_j \sum_{x=0}^{2N} \left( 1 - \frac {x}{N}\right)^n
K_i(x) K_j(x) \frac {\begin{pmatrix}
2N \\
x
\end{pmatrix}}{2^{2N}}.
\]
We can use these expressions to rederive some results given in [K].

We have
\[
({\mathbb P}^n)_{00} = \sum_{x=0}^{2N} \left( 1 - \frac {x}{N}\right)^n \frac
{\begin{pmatrix}
2N \\
x
\end{pmatrix}}{2^{2N}}
\]
and the ``generating function'' for these probabilities, defined by
\[
U(z) \equiv \sum_{n=0}^{\infty} z^n({\mathbb P}^n)_{00}
\]
becomes
\[
U(z) = \sum_{x=0}^{2N} \frac {N}{N(1-z)+xz} \frac {\begin{pmatrix}
2N \\
x
\end{pmatrix}}{2^{2N}}.
\]
In particular $U(1) = \infty$ and then the familiar ``renewal
equation'', see [F] given by
\[
U(z) = F(z)U(z) + 1
\]
where $F(z)$ is the generating function for the probabilities $f_n$
of returning from state $0$ to state $0$ for the first time in $n$
steps
\[
F(z) = \sum_{n=0}^{\infty} z^nf_n
\]
gives
\[
F(z) = 1 - \frac {1}{U(z)}
\]
Therefore we have  $F(1) = 1$, indicating that one returns to state $0$ with
probability one in finite time.

These results allow us to compute the expected time to return to
state $0$.  This expected value is given by $F'(1)$, and we have
\[
F'(z) = \frac {U'(z)}{U^2(z)}.
\]
Since
\[
U'(z) = \sum_{x=0}^{2N} \frac {N(N-x)}{(N(1-z)+xz)^2} \frac {\begin{pmatrix}
2N \\
x
\end{pmatrix}}{2^{2N}}
\]
we get $F'(1) = 2^{2N}$.  The same method shows that any state $i =
0,\dots,2N$ is also recurrent and that the expected time to return to
it is given by
\[
\frac {2^{2N}}{\begin{pmatrix}
2N \\
i
\end{pmatrix}} .
\]
The moral of this story is clear:  if $i=0$ or $2N$, or if $i$ is close to these
values, i.e., we start from a state where most balls are in one urn it
will take on average a huge amount of time to get back to this state.

If on the other hand $i=N$, i.e., we are starting from a very
balanced state, then we will (on average) return to this state fairly
soon. Thus we see how the issues of irreversibility and recurrence are rather
subtle.

\bigskip
In a very precise sense these polynomials are discrete analogs of
those of Hermite in the case of the
real line. For interesting material regarding this section the reader should
consult \cite{A}.

\section{A Chebyshev type example}

The example below illustrates nicely how certain recurrence properties
of the process are related to the presence of point masses in the orthogonality
measure. This is seen by comparing the two integrals at the end of the section.

\bigskip
Consider the matrix
\[
{\mathbb P} = \begin{pmatrix}
0 & 1 & 0 \\
q & 0 & p \\
0 & q & 0 & p \\
& &\ddots &\ddots &\ddots
\end{pmatrix}
\]
with $0 \le p \le 1$ and $q = 1 - p$.  We look for polynomials
$Q_j(x)$ such that
\[
Q_{-1}(x) = 0,\quad Q_0(x) = 1
\]
and if $Q(x)$ denotes the vector
\[
Q(x) \equiv \begin{pmatrix}
Q_0(x) \\
Q_1(x) \\
Q_2(x) \\
\vdots
\end{pmatrix}
\]
we ask that we should have
\[
{\mathbb P}Q(x) = xQ(x).
\]

\bigskip
The matrix ${\mathbb P}$ can be conjugated into a symmetric one and
in this fashion one can
find the explicit expression for these polynomials.

\bigskip
We have
\[
Q_j(x) = \left( \frac {q}{p}\right)^{j/2} \left[(2-2p)T_j \left(
\frac {x}{2\sqrt{pq}}\right) + (2p-1)U_j \left( \frac
{x}{2\sqrt{pq}}\right)\right]
\]
where $T_j$ and $U_j$ are the Chebyshev polynomials of the first and
second kind.

If $p \ge 1/2$ we have
\[
\left( \frac {p}{1-p} \right)^n \int_{-\sqrt{4pq}}^{\sqrt{4pq}}
Q_n(x)Q_m(x) \frac {\sqrt{4pq-x^2}}{1-x^2} dx = \delta_{nm}
\begin{cases}
2(1-p)\pi, &n = 0 \\
2p(1-p)\pi, &n \ge 1
\end{cases}
\]
while if $p \le 1/2$ we get a new phenomenon, namely the presence of
point masses
in the spectral measure

\begin{eqnarray*}
& &\left( \frac {p}{1-p} \right)^n \left[
\int_{-\sqrt{4pq}}^{\sqrt{4pq}} Q_n(x)Q_m(x) \frac
{\sqrt{4pq-x^2}}{1-x^2} dx \right. \\
& &\left. \begin{matrix}
{} \\
{} \\
\end{matrix} + (2-4p)\pi [Q_n(1)Q_m(1) +
Q_n(-1)Q_m(-1)]\right] \\
& &= \delta_{nm} \begin{cases}
2(1-p)\pi, &n = 0 \\
2p(1-p)\pi, &n \ge 1
\end{cases}
\end{eqnarray*}

From a probabilistic point of view these results are very natural.

\section{The Hahn polynomials, Laplace and Bernoulli}

As has been pointed out before, a limitation of this method is given
by the sad fact that
given the matrix ${\mathbb P}$ very seldom can one write down the
corresponding polynomials and their
orthogonality measure. In general there is no reason why "physically
interesting" Markov chains
will give rise to situations where these mathematical objects can be
found explicitly.

\bigskip
The example below shows that one can get lucky: there is a very old
model of the exchange of heat
between two bodies going back to Laplace and Bernoulli, see \cite{F}, page 378.
It turns out that in this case
the corresponding orthogonal polynomials can be determined explicitly.

The Bernoulli--Laplace model for the exchange of heat between two
bodies consists of two urns, labeled 1 and 2.  Initially there are
$W$ white balls in urn~1 and $B$ black balls in urn~2.  The
transition mechanism is as follows:  a ball is picked from each urn and
these two balls are switched.  It is natural to expect that
eventually both urns will have a nice mixture of white and black
balls.

The state of the system at any time is described by $w$, defined to be
the number of white balls in urn~1.
It is clear that we have, for $w = 0,1,\dots,W$
\begin{eqnarray*}
{\mathbb P}_{w,w+1} &= &\frac {W-w}{W} \frac {W-w}{B} \\
{\mathbb P}_{w,w-1} &= &\frac {w}{W} \frac {B-W+w}{B} \\
{\mathbb P}_{w,w} &= &\frac {w}{W} \frac {W-w}{B} + \frac {W-w}{W}
\frac {B-W+w}{B} .
\end{eqnarray*}
Notice that
\[
{\mathbb P}_{w,w-1} + {\mathbb P}_{w,w} + {\mathbb P}_{w,w+1} = 1.
\]

Introduce now the dual Hahn polynomials by means of
\begin{eqnarray*}
R_n(\lambda(x)) &= &{}_3{\tilde F}_2 \left( \left. \begin{matrix}
-n,-x,x-W-B-1 \\
-W,-W
\end{matrix} \right| 1\right) \\
n &= &0,\dots,W;\ x = 0,\dots,W.
\end{eqnarray*}

These polynomials depend in general on one more parameter.

\bigskip

Notice that these are polynomials of degree $n$ in
\[
\lambda(x) \equiv x(x-W-B-1).
\]
One has
\[
{\mathbb P}_{w,w-1}R_{w-1} + {\mathbb P}_{w,w}R_w + {\mathbb
P}_{w,w+1}R_{w+1} = \left( 1 - \frac
{x(B+W-x+1)}{BW}\right) R_w.
\]
This means that for each value of $x = 0,\dots,W$ the vector
\[
\begin{pmatrix}
R_0(\lambda(x)) \\
R_1(\lambda(x)) \\
\vdots \\
R_W(\lambda(x))
\end{pmatrix}
\]
is an eigenvector of the matrix ${\mathbb P}$ with eigenvalue $1 -
\frac {x(B+W-x+1)}{BW}$.  The relevant orthogonality relation is
given by
\[
\pi_j \sum_{x=0}^W R_i(\lambda(x))R_j(\lambda(x))\mu(x) = \delta_{ij}
\]
with
\begin{eqnarray*}
\mu(x) &= &\frac
{w!(-w)_x(-w)_x(2x-W-B-1)}{(-1)^{x+1}x!(-B)_x(x-W-B-1)_{w+1}} \\
\pi_j &= &\frac {(-w)_j}{j!} \frac {(-B)_{w-j}}{(w-j)!}.
\end{eqnarray*}
The Karlin--McGregor representation gives
\[
({\mathbb P}^n)_{ij} = \pi_j \sum_{x=0}^W
R_i(\lambda(x))R_j(\lambda(x))e^n(x)\mu(x)
\]
with $e(x) = 1 - \frac {x(B+W-x+1)}{BW}$.

These results can be used, once again, to get some quantitative
results on this process.
%If we put
%\[
%U^{ij}(z) \equiv \sum_{n=0}^{\infty} ({\mathbb P}^n)_{ij} z^n
%\]
%we get
%\[
%U^{ij}(z) = \pi_j \sum_{x=0}^W R_i(\lambda(x))R_j(\lambda(x)) \frac
%{1}{1-ze(x)} \mu(x)
%\]
%as for $z \approx 1$ the important contribution comes from $x=0$, resulting in
%\[
%\frac {U^1(z)}{U^2(z)} \approx \frac {1}{\pi_j} \frac
%{1}{R_i(\lambda(0))R_j(\lambda(0))\mu(0)}.
%\]

Interestingly enough, these polynomials were considered in great
detail by S.~Karlin and J.~McGregor, \cite{KMcG2} and used by these
authors in the context of a model in genetics
describing fluctuations of gene frequency under the influence of
mutation and selection.
The reader will find useful remarks in \cite{DS}.

\section{The classical orthogonal polynomials and the bispectral problem}

The examples discussed above illustrate the following point: quite
often the orthogonal polynomials
that are associated with important Markov chains belong to the small
class of polynomials
usually referred to as {\em classical}. By this one means that they
satisfy not only three term
recursion relations but that they are also the common eigenfunctions
of some fixed
(usually second order) differential operator.  The search for
polynomials of this kind goes back
at least to \cite{B}. In fact this issue is even older, and was considered
in a paper in 1884 by Routh. See \cite{I,R} for a more complete discussion.

In the context where both variables are
continuous, this problem has been
raised in \cite{DG}. For a view of some related subjects see
\cite{HK}. The reader will
find useful material in \cite{AW,AAR,I}.

\section{Matrix valued orthogonal polynomials}

Here we recall a notion due to M.~G.~Krein, see \cite{K1,K2}.
Given a self adjoint positive
definite matrix valued smooth weight function $W(x)$ with finite
moments, we can consider the skew symmetric bilinear form defined
for any pair of matrix valued polynomial functions $P(x)$ and
$Q(x)$ by the numerical matrix

\[
(P,Q)=(P,Q)_W = \int_{\mathbb R} P(x)W(x)Q^*(x)dx,
\]
where $Q^*(x)$ denotes the conjugate transpose of $Q(x)$. By the
usual construction this leads to the existence of a sequence of
matrix valued orthogonal polynomials with non singular leading
coefficient.

Given an orthogonal sequence $\{P_n (x)\}_{n\ge0}$  one gets by
the usual argument a three term recursion relation
\begin{equation}\label{difference}
xP_n(x) = A_{n}P_{n-1}(x) + B_nP_n(x) + C_{n}P_{n+1}(x),
\end{equation}
where $A_{n}$, $B_{n}$ and $C_{n}$ are matrices and the last one is
nonsingular.

We now turn our attention to an important class of orthogonal
polynomials which we will call {\em classical matrix valued
orthogonal polynomials}.
Very much as in \cite{D,GPT2,GPT3} we say that the weight function is {\em
classical} if there exists a second order ordinary differential
operator $D$ with matrix valued polynomial coefficients $A_j(x)$
of degree less or equal to $j$ of the form
\begin{equation}\label{D}
D=A_2(x)\frac{d^2}{dx^2}+A_1(x)\frac{d}{dx}+A_0(x),
\end{equation}
such that
for an orthogonal
sequence $\{P_n\}$,we have
\begin{equation}\label{autofuncion}
DP_n^*=P_n^*\Lambda_n,
\end{equation}
where $\Lambda_n$ is a real valued matrix.
This form of the eigenvalue equation \eqref{autofuncion} appears
naturally in \cite{GPT1} and differs only superficially with the form
used in \cite{D},
where one uses right handed differential
operators.

During the last few years there has been quite a bit of activity
centered around an effort to produce families of matrix valued
orthogonal polynomials that would satisfy differential equations as
those above.
One of the examples that have resulted from this search, see \cite{G}, will be
particularly useful later on.

\section{Pentadiagonal matrices and matrix valued orthogonal polynomials}

Given a pentadiagonal scalar matrix it is often useful to think of it either in its original
unblocked form or as being made, let us say, of $2 \times 2$ blocks. These two
ways of seeing a matrix, and the fact that matrix operations like
multiplication
can by performed "by blocks" has proved very important in the development of
fast algorithms.

\bigskip
\bigskip

In the case of a birth and death process it is useful to think of a graph like

\vskip2\baselineskip

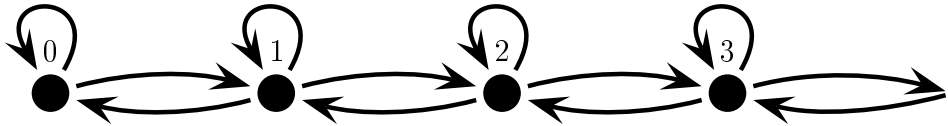
\begin{psmatrix}[rowsep=2.5cm,colsep=3cm]
  \mynode{0} & \mynode{1} & \mynode{2} & \mynode{3} & \pnode{4}
\psset{nodesep=3pt,arcangle=15,labelsep=1ex,linewidth=0.6mm,arrows=->,arrowsize=2mm 3}
\nccurve[angleA=60,angleB=120,ncurv=7]{0}{0}
\nccurve[angleA=60,angleB=120,ncurv=7]{1}{1}
\nccurve[angleA=60,angleB=120,ncurv=7]{2}{2}
%\nccurve[angleA=60,angleB=120,ncurv=7]{4}{4}
%\nccurve[angleA=60,angleB=120,ncurv=7]{6}{6}
\nccurve[angleA=60,angleB=120,ncurv=7]{3}{3}
%\nccurve[angleA=240,angleB=300,ncurv=7]{5}{5}
%\nccurve[angleA=240,angleB=300,ncurv=7]{7}{7}
\ncarc{0}{1} \ncarc{1}{2} \ncarc{2}{3} \ncarc{3}{4}
             \ncarc{1}{0} \ncarc{2}{1} \ncarc{3}{2} \ncarc{4}{3}
\nput{90}{0}{0} \nput{90}{1}{1} \nput{90}{2}{2} \nput{90}{3}{3}
\end{psmatrix}
\vskip2\baselineskip

Suppose that we are dealing with a more complicated Markov chain in
the same probability
space, where the elementary transitions can go beyond "nearest
neighbours". In such a case
the graph may look as follows

\vskip2\baselineskip
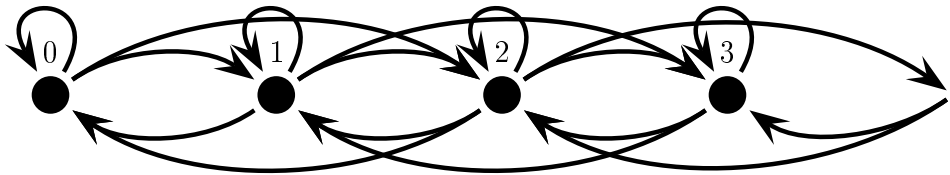
\begin{psmatrix}[rowsep=2.5cm,colsep=3cm]
  \mynode{0} & \mynode{1} & \mynode{2} & \mynode{3} & \pnode{4}
\psset{nodesep=3pt,arcangle=35,labelsep=1ex,linewidth=0.6mm,arrows=->,arrowsize=2mm 3}
\nccurve[angleA=60,angleB=120,ncurv=7]{0}{0}
\nccurve[angleA=60,angleB=120,ncurv=7]{1}{1}
\nccurve[angleA=60,angleB=120,ncurv=7]{2}{2}
%\nccurve[angleA=60,angleB=120,ncurv=7]{4}{4}
%\nccurve[angleA=60,angleB=120,ncurv=7]{6}{6}
\nccurve[angleA=60,angleB=120,ncurv=7]{3}{3}
%\nccurve[angleA=240,angleB=300,ncurv=7]{5}{5}
%\nccurve[angleA=240,angleB=300,ncurv=7]{7}{7}
\ncarc{0}{1} \ncarc{0}{2} \ncarc{1}{2} \ncarc{1}{3} \ncarc{2}{3} \ncarc{2}{4}
             \ncarc{1}{0} \ncarc{2}{0} \ncarc{2}{1} \ncarc{3}{1} \ncarc{3}{2} \ncarc{4}{2} \ncarc{4}{3}
\nput{90}{0}{0} \nput{90}{1}{1} \nput{90}{2}{2} \nput{90}{3}{3}
\end{psmatrix}
\vskip3\baselineskip

The matrix ${\mathbb P}$ going with the graph above is now pentadiagonal. By thinking of it in
the manner mentioned above
we get a block tridiagonal matrix. As an extra bonus, its
off-diagonal blocks are triangular.

\clearpage

The graph below
%%%%%%
\vskip5\baselineskip
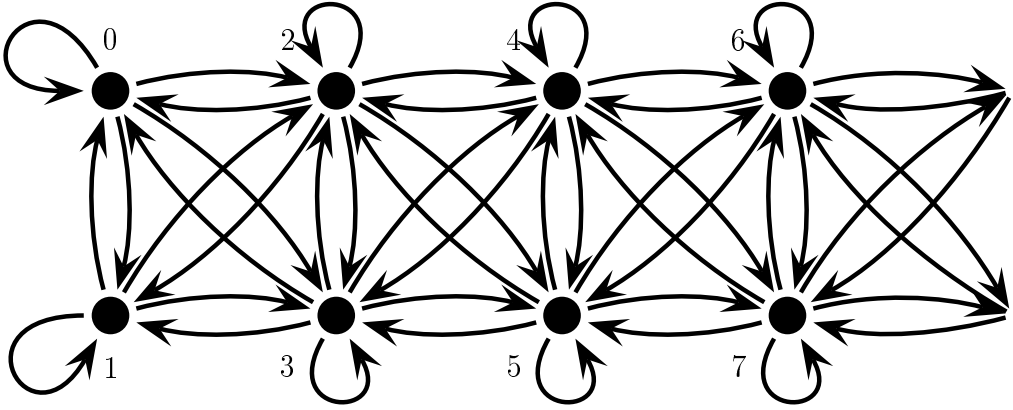
\begin{psmatrix}[rowsep=2.5cm,colsep=3cm]
  \mynode{0} & \mynode{2} & \mynode{4} & \mynode{6} & \pnode{8} \\
  \mynode{1} & \mynode{3} & \mynode{5} & \mynode{7} & \pnode{9}\\
\psset{nodesep=3pt,arcangle=15,labelsep=1ex,linewidth=0.6mm,arrows=->,arrowsize=2mm 3}
\nccurve[angleA=120,angleB=180,ncurv=10]{0}{0}
\nccurve[angleA=180,angleB=240,ncurv=10]{1}{1}
\nccurve[angleA=60,angleB=120,ncurv=7]{2}{2}
\nccurve[angleA=60,angleB=120,ncurv=7]{4}{4}
\nccurve[angleA=60,angleB=120,ncurv=7]{6}{6}
\nccurve[angleA=240,angleB=300,ncurv=7]{3}{3}
\nccurve[angleA=240,angleB=300,ncurv=7]{5}{5}
\nccurve[angleA=240,angleB=300,ncurv=7]{7}{7}
\ncarc{0}{2} \ncarc{2}{4} \ncarc{4}{6} \ncarc{6}{8}
             \ncarc{2}{0} \ncarc{2}{5} \ncarc{4}{2} \ncarc{4}{7} \ncarc{6}{4} \ncarc{8}{6}
\ncarc{0}{1} \ncarc{0}{3} \ncarc{2}{1} \ncarc{4}{3} \ncarc{6}{5} \ncarc{8}{7} \ncarc{6}{9}
\ncarc{1}{0}
\ncarc{1}{3} \ncarc{3}{5} \ncarc{5}{7} \ncarc{7}{9}
             \ncarc{3}{1} \ncarc{5}{3} \ncarc{7}{5}
\ncarc{1}{2} \ncarc{3}{0} \ncarc{3}{4} \ncarc{5}{2} \ncarc{5}{6} \ncarc{7}{4} \ncarc{7}{8}
\ncarc{2}{3} \ncarc{3}{2} \ncarc{4}{5} \ncarc{5}{4} \ncarc{6}{7} \ncarc{7}{6} \ncarc{9}{7} \ncarc{9}{6}
\psset{labelsep=1.7ex}\def\SQRTwo{2.82843}
\nput{90}{0}{0} \nput[labelsep=\SQRTwo ex]{135}{2}{2} \nput[labelsep=\SQRTwo ex]{135}{4}{4} \nput[labelsep=\SQRTwo ex]{135}{6}{6}
\nput{-90}{1}{1} \nput[labelsep=\SQRTwo ex]{-135}{3}{3} \nput[labelsep=\SQRTwo ex]{-135}{5}{5} \nput[labelsep=\SQRTwo ex]{-135}{7}{7}
\end{psmatrix}

\noindent
clearly corresponds to a general
block tridiagonal matrix, with blocks of size $2 \times 2$.

If ${\mathbb P}_{i,j}$ denotes the $i,j$ block of ${\mathbb P}$
we can generate a sequence of $2 \times 2$ matrix valued polynomials
$Q_i(t)$ by imposing
the three term recursion of section 8. By using the notation of
section 2, we would have
\[
{\mathbb P} Q(x) = x Q(x)
\]
where the entries of the column vector $Q(x)$ are now $2 \times 2$ matrices.

Proceeding as in the scalar case, this relation can be iterated to give
\[
{\mathbb P}^n Q(x) = x^n Q(x)
\]
and if we assume the existence of a weight matrix $W(x)$ as in
section 7, with the
property
\[
(Q_j,Q_j)\delta_{i,j} = \int_{\mathbb R} Q_i(x)W(x)Q_j^*(x)dx,
\]
it is then clear that one can get an expression for the $(i,j)$ entry
of the block
matrix ${\mathbb P}^n$ that would look exactly as in the scalar case, namely
%\[
%\begin{pmatrix}
%{\mathbb P}^n_{6,14} & {\mathbb P}_{6,15}^n \\
%{\mathbb P}_{7,14}^n & {\mathbb P}_{7,15}^n
%\end{pmatrix} .
%\]
\[
({\mathbb P}^n)_{ij}  (Q_j,Q_j) = \int x^n Q_i(x)W(x)Q_j^*(x)dx.
\]

Just as in the scalar case, this expression becomes useful when we
can get our hands on the matrix valued
polynomials $Q_i(x)$ and the orthogonality measure $W(x)$. Notice
that we have not discussed conditions on the
matrix ${\mathbb P}$ to give rise to such a measure. For this issue
the reader can consult \cite{DP,D1} and the references in these
papers.

The spectral theory of a scalar double-infinite tridiagonal matrix leads
naturally to a $2 \times 2$ semi-infinite matrix. This has been looked at
in terms of random walks in \cite{P}. In \cite{ILMV} this work is elaborated
further to get a formula that could be massaged to look like the right hand side of the one above.

\section{An explicit example}

Consider the matrix valued polynomials given by the $3$-term recursion relation
\[
A_n\Phi_{n-1}(x) + B_n\Phi_n(x) +  C_n\Phi_{n+1}(x) =
t\Phi_n(x), n \ge 0
\]
\[
\Phi_{-1}(x) = 0,\  \Phi_0(t) = I
\]
where the entries in $A_n$, $B_n$, $ C_n$ are given
below.

Notice that the matrices $A_n$ and $C_n$ are upper and lower
triangular respectively.
\begin{eqnarray*}
A_n^{11} &:= &\frac
{n(\a+n)(\b+2\a+2n+3)}{(\b+\a+2n+1)(\b+\a+2n+2)(\b+2\a+2n+1)} \\
A_n^{12} &:= &\frac {2n(\b+1)}{(\b+2n+1)(\b+\a+2n+2)(\b+2\a+2n+1)} \\
A_n^{21} &:= &0 \\
A_n^{22} &:= &\frac {n(\a+n+1)(\b+2n+3)}{(\b+2n+1)(\b+\a+2n+2)(\b+\a+2n+3)} \\
B_n^{11} &:= &1 + \frac {n(\b+n+1)(\b+2n-1)}{(\b+2n+1)(\b+\a+2n+1)} \\
&- &\frac
{(n+1)(\b+n+2)(\b+2n+1)}{(\b+2n+3)(\b+\a+2n+3)} \\
&- &\frac {2(\b+1)^2}{(\b+2n+1)(\b+2n+3)(\b+2\a+2n+3)} \\
B_n^{12} &:= &\frac {2(\b+1)(\a+\b+n+2)}{(\b+2n+3)(\b+\a+2n+2)(\b+2\a+2n+3)} \\
B_n^{21} &:= &\frac {2(\a+n+1)(\b+1)}{(\b+2n+1)(\b+\a+2n+3)(\b+2\a+2n+3)} \\
B_n^{22} &:= &1 + \frac {n(\b+n+1)(\b+2n+3)}{(\b+2n+1)(\b+\a+2n+2)} \\
&- &\frac {(n+1)(\b+n+2)(\b+2n+5)}{(\b+2n+3)(\b+\a+2n+4)} \\
&+ &\frac {2(\b+1)^2}{(\b+2n+1)(\b+2n+3)(\b+2\a+2n+3)}
\end{eqnarray*}
and finally for the entries of $C_n$ we have
\begin{eqnarray*}
C_n^{11} &:= &\frac
{(\b+n+2)(\b+2n+1)(\b+\a+n+2)}{(\b+2n+3)(\b+\a+2n+2)(\b+\a+2n+3)} \\
C_n^{12} &:= &0 \\
C_n^{21} &:= &\frac {2(\b+1)(\b+n+2)}{(\b+2n+3)(\b+\a+2n+3)(\b+2\a+2n+5)} \\
C_n^{22} &:= &\frac
{(\b+n+2)(\b+\a+n+3)(\b+2\a+2n+3)}{(\b+\a+2n+3)(\b+\a+2n+4)(\b+2\a+2n+5)}.
\end{eqnarray*}

If the matrix $\Psi_0(x)$ is given by
\[
\Psi_0(x) = \begin{pmatrix}
1 & 1 \\
1 & \frac {(\b+2\a+3)x}{\b+1} - \frac {2(\a+1)}{\b+1}
\end{pmatrix}
\]
one can see that the polynomials $\Phi_n(x)$ satisfy the orthogonality relation
\[
\int_0^1 \Phi_i(x)W(x)\Phi_j^*(x)dx = 0 \mbox{ if } i \ne j
\]
where
\[
W(x) = \Psi_0(x) \begin{pmatrix}
(1-x)^{\b}x^{\a+1} & 0 \\
0 & (1-x)^{\b}x^{\a}
\end{pmatrix} \Psi_0^*(x).
\]
The polynomials $\Phi_n(x)$ are ``classical'' in the sense that they are
eigenfunctions of a fixed second order differential operator.  More
precisely, we
have
\[
{\mathcal F}\Phi_n^* = \Phi_n^*\Lambda_n
\]
where
\[
\Lambda_n = \begin{pmatrix}
-n^2-(\a+\b+2)n+\a+1+\frac {\b+1}{2} & 0 \\
0 & -n^2-(\a+\b+3)n
\end{pmatrix}
\]
and ${\mathcal F}$ is given by
\begin{eqnarray*}
{\mathcal F} &= &x(1-x)\left( \frac {d}{dx} \right)^2 \\
&+ &\begin{pmatrix}
\frac {(\a+1)(\b+2\a+5)}{\b+2\a+3} - (\a+\b+3)x & \frac
{2\a+2}{2\a+\b+3} + x \\
\frac {\b+1}{\b+2\a+3} & \frac {(\a+2)\b+2\a^2+5\a+4}{\b+2\a+3} - (\a+\b+4)x
\end{pmatrix} \frac {d}{dx} \\
&+ &\begin{pmatrix}
\a+1+\frac {\b+1}{2} & 0 \\
0 & 0
\end{pmatrix} I.
\end{eqnarray*}

As mentioned earlier this is the reason why this example has surfaced
recently, see \cite{G}.
An explicit expression for the polynomials thenselves is given in
\cite{T2}, corollary 3.

\bigskip
Now we make two observations

\begin{itemize}
\item[1)] the entries of the corresponding pentadigonal matrix are
all non-negative.
\item[2)] the sum of the entries on any given row are all equal to $1$.
\end{itemize}

\bigskip
\bigskip
\bigskip

This allows for an immediate probabilistic interpretation of the pentadiagonal
matrix as the one step transition probability matrix for a Markov chain
whose state space could be visualized by the graph given below.

%%%%%%
\vskip3\baselineskip
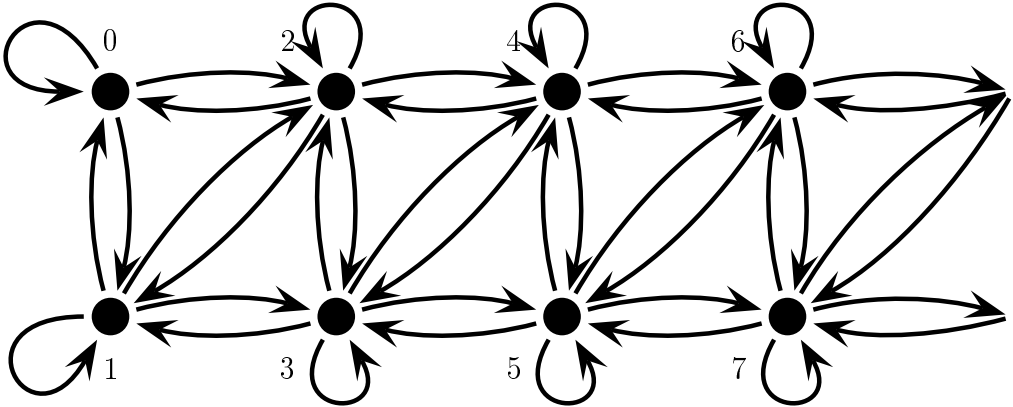
\begin{psmatrix}[rowsep=2.5cm,colsep=3cm]
  \mynode{0} & \mynode{2} & \mynode{4} & \mynode{6} & \pnode{8} \\
  \mynode{1} & \mynode{3} & \mynode{5} & \mynode{7} & \pnode{9}\\
\psset{nodesep=3pt,arcangle=15,labelsep=1ex,linewidth=0.6mm,arrows=->,arrowsize=2mm 3}
\nccurve[angleA=120,angleB=180,ncurv=10]{0}{0}
\nccurve[angleA=180,angleB=240,ncurv=10]{1}{1}
\nccurve[angleA=60,angleB=120,ncurv=7]{2}{2}
\nccurve[angleA=60,angleB=120,ncurv=7]{4}{4}
\nccurve[angleA=60,angleB=120,ncurv=7]{6}{6}
\nccurve[angleA=240,angleB=300,ncurv=7]{3}{3}
\nccurve[angleA=240,angleB=300,ncurv=7]{5}{5}
\nccurve[angleA=240,angleB=300,ncurv=7]{7}{7}
\ncarc{0}{2} \ncarc{2}{4} \ncarc{4}{6} \ncarc{6}{8}
             \ncarc{2}{0} \ncarc{4}{2} \ncarc{6}{4}
\ncarc{0}{1} \ncarc{2}{1} \ncarc{4}{3} \ncarc{6}{5} \ncarc{8}{6}
\ncarc{1}{0}
\ncarc{1}{3} \ncarc{3}{5} \ncarc{5}{7} \ncarc{7}{9}
             \ncarc{3}{1} \ncarc{5}{3} \ncarc{7}{5}
\ncarc{1}{2} \ncarc{3}{4} \ncarc{5}{6} \ncarc{7}{8}
\ncarc{2}{3} \ncarc{3}{2} \ncarc{4}{5} \ncarc{5}{4} \ncarc{6}{7} \ncarc{7}{6} \ncarc{8}{7} \ncarc{9}{7}
\psset{labelsep=1.7ex}\def\SQRTwo{2.82843}
\nput{90}{0}{0} \nput[labelsep=\SQRTwo ex]{135}{2}{2} \nput[labelsep=\SQRTwo ex]{135}{4}{4} \nput[labelsep=\SQRTwo ex]{135}{6}{6}
\nput{-90}{1}{1} \nput[labelsep=\SQRTwo ex]{-135}{3}{3} \nput[labelsep=\SQRTwo ex]{-135}{5}{5} \nput[labelsep=\SQRTwo ex]{-135}{7}{7}
\end{psmatrix}

I find it rather remarkable that this example which was produced for an
entirely different purpose should have this extra property. Finding an
appropriate combinatorial mechanism, maybe in terms of urns, that goes
along with this example remains an interesting challenge.

Two final observations dealing with these state spaces that can be analyzed
using matrix valued orthogonal polynomials. If we were using matrix valued
polynomials of size $N$ we would have as state space a semiinfinite network
consisting of $N$ (instead of two) parallel collection of non-negative integers
with connections going from each of the $N$ states on each vertical rung to every one in the same rung and the two neighbouring ones. The examples in \cite{GPT1}
give instances of this situation with a rather local connection pattern.

In the case of $N=2$ one could be tempted to paraphrase a well known paper and
say that "it has not escaped our notice that" some of these models could be used to study transport phenomena along a DNA segment.

\section{Another example}

Here we consider a different example of matrix valued orthogonal
polynomials whose
block tridiagonal matrix can be seen as a scalar pentadiagonal matrix
with non-negative
elements. In this case the sum of the elements in the rows of this
scalar matrix is not
identically one, but this poses no problem in terms of a
Karlin--McGregor type representation
formula for the entries of the powers ${\mathbb P}^n$.

This example has the important property that the orthogonality weight
matrix $W(x)$, as well as the polynomials themselves are
explicitly known. This is again, a classical situation, see \cite{CG}.

Consider the block tridiagonal matrix
\[
\begin{pmatrix}
B_0 & I \\
A_1 & B_1 & I \\
0 & A_2 & B_2 & I \\
& &\ddots &\ddots &\ddots
\end{pmatrix}
\]
with $2 \times 2$ blocks given as follows
\begin{eqnarray*}
B_0 &= &\frac {1}{2} I,\ B_n = 0 \mbox{ if } n \ge 1 \\
A_n &= &\frac {1}{4} I \mbox{ if } n \ge 1.
\end{eqnarray*}
In this case one can compute explicitly the matrix valued polynomials
$P_n$ given by
\[
A_nP_{n-1}(x) + B_nP_n(x) + P_{n+1}(x) = xP_n(x), \ P_{-1}(x) = 0,\ P_0(x) = I.
\]
One gets
\[
P_n(x) = \frac {1}{2^n} \begin{pmatrix}
U_n(x) & -U_{n-1}(x) \\
-U_{n-1}(x) & U_n(x)
\end{pmatrix}
\]
where $U_n(x)$ are the Chebyshev polynomials of the second kind.

\bigskip
\bigskip

The orthogonality measure is read off from the identity
\[
\frac {4^i}{\pi} \int_{-1}^1 P_i(x) \frac {1}{\sqrt{1-x^2}} \begin{pmatrix}
1 & x \\
x & 1
\end{pmatrix} P_j(x)dx = \delta_{ij} I.
\]
We get, for $n = 0,1,2,\dots$
\[
\frac {4^i}{\pi} \int_{-1}^1 x^nP_i(x) \frac {1}{\sqrt{1-x^2}} \begin{pmatrix}
1 & x \\
x & 1
\end{pmatrix} P_j(x)dx = ({\mathbb P}^n)_{ij}
\]
where, as above,
\[
({\mathbb P}^n)_{ij}
\]
stands for the $i,j$ block of the matrix ${\mathbb P}^n$.

In this way one can compute the entries of the powers ${\mathbb P}^n$
with ${\mathbb P}^n$ thought of as a pentadiagonal matrix, namely

\bigskip
\bigskip

\[
{\mathbb P} = \begin{pmatrix}
0 & \frac {1}{2} & 1 & 0 & 0 & 0 \\
\frac {1}{2} & 0 & 0 & 1 & 0 & 0 & \ddots \\
\frac {1}{4} & 0 & 0 & 0 & 1 & 0 & \ddots \\
0 & \frac {1}{4} & 0 & 0 & 0 & 1 & \ddots \\
& 0 & \frac {1}{4} & 0 & 0 & 0 & \ddots \\
& & 0 & \frac {1}{4} & 0 & 0 & \ddots \\
& & & \ddots & \ddots & \ddots & \ddots
\end{pmatrix}.
\]

\bigskip
\bigskip
\bigskip

This example goes along with the following graph.

%%%%%
\vskip3\baselineskip
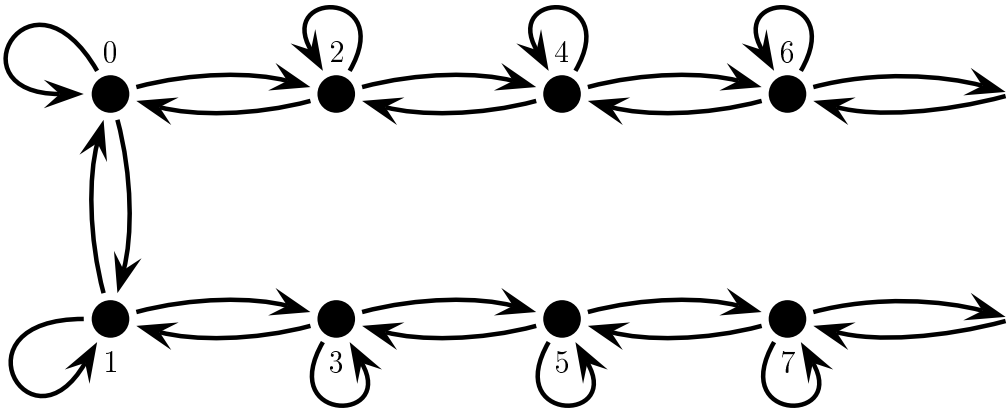
\begin{psmatrix}[rowsep=2.5cm,colsep=3cm]
  \mynode{0} & \mynode{2} & \mynode{4} & \mynode{6} & \pnode{8} \\
  \mynode{1} & \mynode{3} & \mynode{5} & \mynode{7} & \pnode{9}\\
\psset{nodesep=3pt,arcangle=15,labelsep=1ex,linewidth=0.6mm,arrows=->,arrowsize=2mm 3}
\nccurve[angleA=120,angleB=180,ncurv=10]{0}{0}
\nccurve[angleA=180,angleB=240,ncurv=10]{1}{1}
\nccurve[angleA=60,angleB=120,ncurv=7]{2}{2}
\nccurve[angleA=60,angleB=120,ncurv=7]{4}{4}
\nccurve[angleA=60,angleB=120,ncurv=7]{6}{6}
\nccurve[angleA=240,angleB=300,ncurv=7]{3}{3}
\nccurve[angleA=240,angleB=300,ncurv=7]{5}{5}
\nccurve[angleA=240,angleB=300,ncurv=7]{7}{7}
\ncarc{0}{2} \ncarc{2}{4} \ncarc{4}{6} \ncarc{6}{8}
             \ncarc{2}{0} \ncarc{4}{2} \ncarc{6}{4} \ncarc{8}{6}
\ncarc{0}{1}
\ncarc{1}{0}
\ncarc{1}{3} \ncarc{3}{5} \ncarc{5}{7} \ncarc{7}{9}
             \ncarc{3}{1} \ncarc{5}{3} \ncarc{7}{5} \ncarc{9}{7}
\nput{90}{0}{0} \nput{90}{2}{2} \nput{90}{4}{4} \nput{90}{6}{6}
\nput{-90}{1}{1} \nput{-90}{3}{3} \nput{-90}{5}{5} \nput{-90}{7}{7}
\end{psmatrix}

\section{A few more challenges}

We have already pointed out a few challenges raised by our attempt to extend
the Karlin-McGregor representation beyond its original setup. Here we
list a few
more open problems. The reader will undoubtedly come up with many more.

\bigskip
Is there a natural version of the models introduced by
Bernoulli-Laplace and by P. and T. Ehrenfest whose solution features matrix valued polynomials?

\bigskip
Is it possible to modify the simplest Chebyshev type examples in \cite{D1}
to accomodate cases where some of the blocks in the tridiagonal matrix
give either absorption or reflection boundary conditions?

\bigskip
One could consider the emerging class of polynomials of several variables
and find here interesting instances where the state space is higher dimensional.
For a systematic study of polynomials in several variables one should consult \cite{DX} as well as
the work on Macdonald polynomials of various kinds, see \cite{M}. A look at the pioneering work of Tom Koornwinder, see for instance \cite{Ko}, is always a very good idea.

\bigskip
After this paper was finished I came up with two independent sources of multivariable polynomials of the type alluded to in the previous paragraph. One is the series of papers by Hoare and Rahman, see \cite{HR1,HR2,HR3,HR4}. The other one deals with the papers \cite{Mi} and \cite{IX,GI}.

\bigskip
In queueing theory one finds the notion of Quasi-Birth-and-Death processes, see \cite{LR1,N}. Within those that are nonhomogeneous one could find examples where
the general approach advocated here might be useful.

\end{document}